\documentclass[a4paper,10pt]{article}

\usepackage{amsmath}
\usepackage{amssymb}
\usepackage{amsthm}
\usepackage{mathrsfs}
\usepackage{color}
\usepackage[dvips]{graphicx}
\usepackage[all]{xy}

\theoremstyle{plain}
\newtheorem{them}{Theorem}[section]
\newtheorem{lemma}[them]{Lemma}
\newtheorem{prop}[them]{Proposition}

\theoremstyle{definition}
\newtheorem{defi}[them]{Definition}
\newtheorem{exam}[them]{Example}
\newtheorem{rema}[them]{Remark}

\newtheorem*{mthem}{Main Theorem}

\newcommand{\lmap}[3]{#1:#2 \longrightarrow #3}
\newcommand{\map}[3]{#1:#2 \rightarrow #3}

\begin{document}

\title{Coverings of small categories and nerves}
\author{Kazunori Noguchi \thanks{noguchi@math.shinshu-u.ac.jp}}
\date{}
\maketitle
\begin{abstract}
We prove a certain proposition which states a relationship between coverings of small categories and nerves. As its application, we prove that for a covering $\map{P}{E}{B}$ of finite categories, the zeta function of $E$ is the zeta function of $B$ to the number of sheet of $P$. Moreover, we prove the formula $\chi(E)=\chi(F)\chi(B)$ for Euler characteristic of categories and coverings.
\end{abstract}


\thispagestyle{empty}

\section{Introduction}

A covering space is very useful and interesting tools in geometry. For instance, it is used for calculations of fundamental groups and it has an analogy of Galois theory \cite{Hat02}. The notion of coverings is also defined for small categories. Many people have studied about it, for example \cite{BH99}, \cite{CM} and \cite{Tan}. In this paper, we show the following proposition.
\begin{prop}\label{fiber0}
Let $\map{P}{E}{B}$ be a covering of small categories and let $b$ be an object of $B$. For any $n\ge 0$, $N_n(E)$ is bijective to $\coprod_{x\in P^{-1}(b)} N_n(B)$
where $$N_n(E)=\{\xymatrix{(x_0\ar[r]^{f_1}&x_1\ar[r]^{f_2}&\dots\ar[r]^{f_n}&x_n)} \text{ in } E \}$$ and $N_n(B)$ is also defined in the same way.
\end{prop}
As its application, we prove that for a covering $\map{P}{E}{B}$ of finite categories, the zeta function of $E$ is the zeta function of $B$ to the number of sheet of $P$. Moreover, we prove the formula $\chi(E)=\chi(F)\chi(B)$ for Euler characteristic of categories and coverings.

First, we give a historical background of the zeta function of a finite category and coverings of finite categories.
In \cite{NogA}, the zeta function of a finite category $I$ was defined by $$\zeta_I(z)=\exp\left( \sum_{n=1}^{\infty} \frac{\# N_n(I)}{n} z^n\right)$$
and it was shown that the zeta function of a connected finite groupoid is a rational function and for a covering of finite groupoids $\map{P}{E}{B}$, the inverse zeta function of $B$ divides the inverse zeta function of $E$ (Proposition 2.3 and Proposition 2.7 of \cite{NogA}). 

In this paper, we generalize Proposition 2.7 of \cite{NogA} as follows.
\begin{mthem}Suppose $\map{P}{E}{B}$ is a covering of finite categories and let $b$ be an object of $B$. Then, the zeta function of $E$ is the zeta function of $B$ to the number of sheet of $P$, that is,
$$\zeta_E(z)=\big(\zeta_B(z)\big)^{\# P^{-1}(b)}.$$
\end{mthem}Note that the number $\#P^{-1}(b)$ does not depend on the choice of $b$ (Proposition \ref{sheet}). It is an analogue of Corollary 1 of  \S 2 of \cite{ST96}.

Next, we recall Euler characteristic of categories.

The Euler characteristic of a finite category $\chi_L$ was defined by Leinster \cite{Leia}. After this paper, various definitions of Euler characteristics for categories were defined, series Euler characteristic $\chi_{\sum}$ by Berger-Leinster \cite{Leib}, $L^2$-Euler characteristic $\chi^{(2)}$ by Fiore-L\"uck-Sauer \cite{FLS},  extended $L^2$-Euler characteristic $\chi^{(2)}_{\mathrm{ex}}$ \cite{Nog} and Euler characteristic of $\mathbb{N}$-filtered acyclic categories $\chi_{\mathrm{fil}}$ \cite{Nog11} by the author . For a finite acyclic category, four Euler characteristics $\chi_L,\chi_{\sum},\chi^{(2)}$ and $\chi_{\mathrm{fil}}$ coincide (see, for example, Introduction of \cite{Nog}). A small category is \textit{acyclic} if every endomorphism and every isomorphism is an identity morphism.  For a finite groupoid $G$, three Euler characteristic $\chi_L,\chi_{\sum}$ and $\chi^{(2)}$ coincide and the value is $$\sum_{x\in \mathrm{Ob}(G)/\cong} \frac{1}{\#\mathrm{Aut}(x)}$$
where this sum runs over all isomorphism classes of objects of $G$ (Example 2.7 of \cite{Leia}, Theorem 3.2 of \cite{Leib} and Example 5.12 of \cite{FLS}). 

A topological fibration $F \hookrightarrow E\rightarrow B$ under certain suitable hypothesis satisfies the equation $$\chi(E)=\chi(F)\chi(B).$$ An analogue of this formula for Euler characteristic of categories  and categorical fibrations was considered in \cite{Leia} and \cite{FLS}. In \cite{Leia}, such formula was found for the Grothendieck construction (Proposition 2.8 of \cite{Leia}). In \cite{FLS}, such formula for coverings of connected finite groupoids and isofibrations of connected finite groupoids were found (Theorem 5.30 and Theorem 5.37 of \cite{FLS}).

In this paper, we consider such formula for $\chi_{\sum}$ and $\chi_{\mathrm{fil}}$ and coverings.
\begin{mthem}

\begin{enumerate}
\item Let $\map{P}{E}{B}$ be a covering of finite categories and let $b$ be an object of $B$. Then, $E$ has series Euler characteristic if and only if $B$ has series Euler characteristic. In this case, we have
$$\chi_{\sum}(E)=\chi_{\sum}(P^{-1}(b)) \chi_{\sum}(B).$$

\item Suppose $(A,\mu_A)$ and $(B,\mu_B)$ are $\mathbb{N}$-filtered acyclic categories and $b$ is an object of $B$ and $\map{P}{A}{B}$ be a covering whose fiber is finite satisfying the equation $\mu_A(x)=\mu_B(P(x))$ for any object $x$ of $A$ . Then, $(A,\mu_A)$ has Euler characteristic $\chi_{\mathrm{fil}}(A,\mu_A)$ if and only if $B$ has Euler characteristic $\chi_{\mathrm{fil}}(B,\mu_B)$. In this case, we have
$$\chi_{\mathrm{fil}}(A,\mu_A)=\chi_{\mathrm{fil}}(P^{-1}(b),\mu) \chi_{\mathrm{fil}}(B,\mu_B)$$
for any $\mathbb{N}$-filtration $\mu$ of $P^{-1}(b)$.
\end{enumerate}
\end{mthem}

This paper is organized as follows.

In section \ref{main}, we investigate relationships between coverings of small categories and nerves.

In section \ref{application}, we prove our main theorem.

In section \ref{example}, we give some examples of coverings of small categories.

\section{Coverings and Nerves}\label{main}

In this section, we investigate relationships between coverings of small categories and nerves.

Here, let us recall a covering of small categories \cite{BH99}.

Let $C$ be a small category. Then, $C$ is \textit{connected} if there exists a zig-zag sequence of morphisms in $C$
$$\xymatrix{x\ar[r]^{f_1}&x_1&x_2\ar[r]^{f_3}\ar[l]_{f_2}&\dots&y\ar[l]_(0.4){f_n}}$$
for any objects $x$ and $y$ of $C$. We do not have to care about the direction of the last morphism $f_n$ since we can insert an identity morphism to the sequence. For an object $x$ of $C$, let $S(x)$ is the set of morphisms of $C$ whose source is $x$
$$S(x)=\{\map{f}{x}{*}\in \mathrm{Mor}(C)\}$$
and $T(x)$ is the set of morphisms of $C$ whose target is $x$
$$T(x)=\{\map{g}{*}{x}\in \mathrm{Mor}(C)\}.$$ Suppose $E$ and $B$ are small categories and $B$ is connected. Then, a functor $\map{P}{E}{B}$ is a \textit{covering} if the following two restrictions of $P$
$$P:S(x)\longrightarrow S(P(x))$$ $$P:T(x)\longrightarrow T(P(x))$$
are bijections for any object $x$ of $E$. This condition is an analogue of the condition of an unramified covering of graphs (see \cite{ST96}). For an object $b$ of $B$, the inverse image $P^{-1}(b)$ of the restriction of $P$ with respect to objects
$$P^{-1}(b)=\{x\in \mathrm{Ob}(E)\mid P(x)=b\}$$ is called the \textit{fiber} of $b$. The cardinality of $P^{-1}(b)$ is called \textit{the number of sheet of }$P$ and it does not depend on the choice of $b$ since the base category $B$ is connected (see Proposition \ref{sheet}). In particular, a covering of groupoids was studied in \cite{May99}. Applying the classifying space functor $B$ to a covering $\map{P}{E}{B}$, we have the covering space $BP$ in the topological sense (see \cite{Tan}) and examples are given in Section \ref{example}.

The following proposition is briefly introduced in \cite{BH99} with no proof, but this proposition is very important in this paper. So we give a proof of it to make sure.

\begin{prop}\label{sheet}
Let $\map{P}{E}{B}$ be a covering of small categories. Then, for any objects $b$ and $b'$ of $B$, $P^{-1}(b)$ is bijective to $P^{-1}(b').$
\begin{proof}
It suffices to show that $P^{-1}(b)$ is bijective to $P^{-1}(b')$ if there exists a morphism $f:b\rightarrow b'$. Indeed, if it is proven, then we have for any objects $b$ and $b'$, we have a zig-zag sequence 
$$\xymatrix{b\ar[r]&b_1&b_2\ar[r]\ar[l]&\dots&b'\ar[l]},$$
so that we obtain 
$$P^{-1}(b)\cong P^{-1}(b_1)\cong\cdots \cong P^{-1}(b').$$

Suppose there exists a morphism $f:b\rightarrow b'$. For any $x$ of $P^{-1}(b)$, we have the bijection $\map{P}{S(x)}{S(b)}$. Since $f$ belongs to $S(b)$, there exists a unique morphism $\map{f_x}{x}{y_x}$ of $S(x)$ such that $P(f_x)=f$. Define a map $\map{F_f}{P^{-1}(b)}{P^{-1}(b')}$ by $F_f(x)=y_x$. Then, $F_f$ is a bijection. We first show it is injective. If $x\not = x'$, then $y_x\not =y_{x'}$. Indeed, then if we assume $y_x=y_{x'}$
$$\xymatrix{x\ar[dr]^{f_x}&&x' \ar[ld]_{f_{x'}}\\&y_x=y_{x'},}$$
$f_x$ and $f_{x'}$ belong to $T(y_x)$ and they are different morphisms, but the map $\map{P}{T(y_x)}{T(b')}$ is bijective and $P(f_x)=P(f_{x'})=f$. This contradiction implies $y_x \not =y_{x'}$. Next, we show $F_f$ is a surjection. For any $y$ of $P^{-1}(b')$, there exists a unique morphism $\map{g}{x}{y}$ such that $P(g)=f$ since the map $\map{P}{T(y)}{T(b')}$ is bijective. Then, $P(x)=b$, so $x$ belongs to $P^{-1}(b)$. Hence, $F_f(x)=y$.
\end{proof}
\end{prop}

The following lemma is clear, but this formulation will make the proof of Proposition \ref{target} easier to understand.

\begin{lemma}\label{decompose}
Suppose $\map{f}{X}{Y}$ is a bijection and $X=\coprod_{\lambda \in \Lambda} A_{\lambda}$ and $Y=\coprod_{\lambda \in \Lambda} B_{\lambda}$ and for each restriction $f|_{A_{\lambda}}$, its image belongs to $B_{\lambda}$. Then, each restriction $f|_{A_{\lambda}}$ is a bijection.
$$\xymatrix{X\ar[d]_f^{\cong}&=&\coprod_{\lambda \in \Lambda} A_{\lambda}\ar[d]_f^{\cong}&\supset&A_{\lambda} \ar[d]_{f|_{A_{\lambda}}}^{\cong}\\Y&=&\coprod_{\lambda \in \Lambda} B_{\lambda}&\supset&B_{\lambda}  }$$
\end{lemma}

Let $C$ be a small category and $x$ be an object of $C$. Then, let $N_n(C)_x$ be the set of chains of morphisms in $C$ of length $n$ whose target is $x$
$$N_n(C)_x=\{\xymatrix{(x_0\ar[r]^{f_1}&x_1\ar[r]^{f_2}&\dots\ar[r]^{f_n}&x_n)} \text{ in } C \mid x_n=x\}.$$

\begin{prop}\label{target}
Let $\map{P}{E}{B}$ be a covering of small categories. Then, for any object $b$ of $B$ and any $x$ of $P^{-1}(b)$ and $n\ge 0$, $N_n(E)_x$ is bijective to $N_n(B)_b$.
\begin{proof}
We prove this proposition by induction on $n$. If $n=0$, we have
$$N_0(E)_x=\{1_x \}\cong \{1_b\}=N_0(B)_b.$$

Suppose it is true for $n$. Then we have 
\begin{eqnarray}
N_{n+1}(E)_x&=&\coprod_{y\in \mathrm{Ob}(E)} N_n(E)_y\times \mathrm{Hom}_E(y,x) \nonumber \\
&=&\coprod_{b_i\in \mathrm{Ob}(B)}  \coprod_{y_i\in P^{-1}(b_i)}N_n(E)_{y_i}\times \mathrm{Hom}_E(y_i,x) \nonumber  \\
&\cong &\coprod_{b_i\in \mathrm{Ob}(B)}  \coprod_{y_i\in P^{-1}(b_i)}N_n(B)_{b_i}\times \mathrm{Hom}_E(y_i,x) \nonumber  \\
&\cong &\coprod_{b_i\in \mathrm{Ob}(B)}  N_n(B)_{b_i}\times\bigg( \coprod_{y_i\in P^{-1}(b_i)} \mathrm{Hom}_E(y_i,x) \bigg). \label{siki1}
\end{eqnarray}
We have the following diagram $$\xymatrix{T(x)=\displaystyle  \coprod_{b_i \in \mathrm{Ob}(B)}\bigg( \coprod_{y_i\in P^{-1}(b_i)} \mathrm{Hom}_E(y_i,x) \bigg) \ar[d]_P^{\cong}&\displaystyle   \coprod_{y_i\in P^{-1}(b_i)} \mathrm{Hom}_E(y_i,x)  \ar[d]_P \ar@{_{(}->}[l] \\
T(b)= \displaystyle  \coprod_{b_i \in \mathrm{Ob}(B)} \mathrm{Hom_B}(b_i,b)&\mathrm{Hom_B}(b_i,b) \ar@{_{(}->}[l] .} $$
Lemma \ref{decompose} implies 
$$\lmap{P}{\coprod_{y_i\in P^{-1}(b_i)} \mathrm{Hom}_E(y_i,x)}{\mathrm{Hom}_B(b_i,b)}$$
is a bijection since for each $\map{f}{y_i}{x}$, $\map{P(f)}{b_i}{b}$ belongs to $\mathrm{Hom}_B(b_i,b)$.
Hence, the equation \eqref{siki1} is
\begin{eqnarray*}
\eqref{siki1}&\cong&\coprod_{b_i\in \mathrm{Ob}(B)} N_n(B)_{b_i}\times \mathrm{Hom}_B(b_i,b) \\
&=&N_{n+1}(B)_b.
\end{eqnarray*}
\end{proof}
\end{prop}

\begin{prop}\label{fiber}
Let $\map{P}{E}{B}$ be a covering of small categories and let $b$ be an object of $B$. For any $n\ge 0$, $N_n(E)$ is bijective to $\coprod_{x\in P^{-1}(b)} N_n(B)$.
\begin{proof}
When $n=0$, Proposition \ref{sheet} implies
\begin{eqnarray*}
N_0(E)&=&\coprod_{b\in \mathrm{Ob}(B)} P^{-1}(b) \\
&\cong& P^{-1}(b)\times N_0(B) \\
&\cong & \coprod_{x\in P^{-1}(b)} N_0(B).
\end{eqnarray*}

For $n\ge 1$, Proposition \ref{target} implies
\begin{eqnarray}
N_n(E)&=&\coprod_{x\in \mathrm{Ob}(E)} N_{n-1}(E)_x\times S(x) \nonumber \\
&\cong & \coprod_{b\in \mathrm{Ob}(B)} \coprod_{x\in P^{-1}(b)} N_{n-1}(B)_b\times S(b) \label{siki2}
\end{eqnarray}
Since the number of sheet of $P$ does not depend on the choice of $b$ (Proposition \ref{sheet}), we have \eqref{siki2} is
\begin{eqnarray*}
\eqref{siki2}&\cong & \coprod_{x\in P^{-1}(b)}  \coprod_{b\in \mathrm{Ob}(B)} N_{n-1}(B)_b\times S(b) \\
&\cong &\coprod_{x\in P^{-1}(b)} N_n(B).
\end{eqnarray*}
\end{proof}
\end{prop}

The two propositions above hold when nerves are \textit{non-degenerate}, it means that we do not use identity morphisms. Let $C$ be a small category and let $x$ and $y$ be objects of $C$. We define the following symbols by
$$\overline{S(x)}=S(x)\setminus\{1_x\},\overline{T(x)}=T(x)\setminus\{1_x\}$$
$$\overline{\mathrm{Hom}_C}(x,y)=\begin{cases} \mathrm{Hom}_C(x,y)\setminus\{1_x\} & \text{if}\ x=y \\ \mathrm{Hom}_C(x,y) & \text{if} \  x\not =y  \end{cases}$$
$$\overline{N_n}(C)=\{\xymatrix{(x_0\ar[r]^{f_1}&x_1\ar[r]^{f_2}&\dots\ar[r]^{f_n}&x_n)} \text{ in } C\mid f_i\not=1 \}$$
$$\overline{N_n}(C)_x=\{\xymatrix{(x_0\ar[r]^{f_1}&x_1\ar[r]^{f_2}&\dots\ar[r]^{f_n}&x_n)} \text{ in } C\mid f_i\not=1, x_n=x \}.$$
Note that $N_0(C)=\overline{N_0}(C)$.

\begin{prop}\label{n-target}
Let $\map{P}{E}{B}$ be a covering of small categories. Then, for any object $b$ of $B$ and any $x$ of $P^{-1}(b)$ and $n\ge 0$, $\overline{N_n}(E)_x$ is bijective to $\overline{N_n}(B)_b$.
\begin{proof}
If we replace the symbols in the proof of Proposition \ref{target} by the symbols with bars above, we can use the same proof.
\end{proof}
\end{prop}

\begin{prop}\label{n-fiber}
Let $\map{P}{E}{B}$ be a covering of small categories and $b$ be an object of $B$. For any $n\ge 0$, $\overline{N_n}(E)$ is bijective to $\coprod_{x\in P^{-1}(b)} \overline{N_n}(B)$.
\end{prop}

\section{Applications}\label{application}

In this section, we prove our main theorem.

The following is an analogue of Corollary 1 of  \S 2 of \cite{ST96}.
\begin{them}\label{covering}
Let $\map{P}{E}{B}$ be a covering of finite categories and  let $b$ be an object of $B$. Then, we have $$\zeta_E(z)=\big(\zeta_B(z)\big)^{\# P^{-1}(b)}.$$
\begin{proof}
The definition of the zeta function of a finite category and Proposition \ref{fiber} directly imply this fact, that is,
\begin{eqnarray*}
\zeta_E(z)&=&\exp\bigg(\sum_{n=1}^{\infty} \frac{\# N_n(E)}{n} z^n\bigg) \\
&=&\exp\bigg(\sum_{n=1}^{\infty} \frac{\#P^{-1}(b) \# N_n(B)}{n} z^n\bigg)  \\
&=&\big(\zeta_B(z)\big)^{\# P^{-1}(b)}.
\end{eqnarray*}
\end{proof}
\end{them}

Let $I$ be a finite category. Then, $I$ has \textit{series Euler characteristic} \cite{Leib} if and only if the rational function
$$f_I(t)=\frac{\text{sum}(\text{adj}(E-(A_I-E)t ))}{\text{det}(E-(A_I-E)t )}$$
can be substituted $-1$ to $t$ where $A_I$ is an $n\times n$-matrix, called \textit{adjacency matrix}, whose $(i,j)$-entry is the number of morphisms from $x_i$ to $x_j$ when 
$$\mathrm{Ob}(I)=\{x_1,x_2,\dots ,x_n\}.$$ If $I$ has series Euler characteristic , then the series Euler characteristic $\chi_{\sum}(I)$ of $I$ is defined by $f_I(-1)$. This rational function is the rational expression of the power series $\sum^{\infty}_{n=0} \# \overline{N_n}(I) t^n$ (Theorem 2.2 of \cite{Leib}).

A \textit{discrete category} consists of only objects and identity morphisms. For a covering $\map{P}{E}{B}$, its fiber is a discrete category when we regard it as a category.

\begin{them}
Let $\map{P}{E}{B}$ be a covering of finite categories and let $b$ be an object of $B$. Then, $E$ has series Euler characteristic if and only if $B$ has series Euler characteristic. In this case, we have
$$\chi_{\sum}(E)=\chi_{\sum}(P^{-1}(b)) \chi_{\sum}(B).$$
\begin{proof}
We give two types of proofs. 

The first one is a proof by Proposition \ref{n-fiber}. Proposition \ref{n-fiber}  and Theorem 2.2 of \cite{Leib} imply
\begin{eqnarray*}
\sum_{n=0}^{\infty} \# \overline{N_n}(E) t^n&=&\# P^{-1}(b) \sum_{n=0}^{\infty} \# \overline{N_n}(B) t^n \\
&=&\# P^{-1}(b) \frac{\text{sum}(\text{adj}(E-(A_B-E)t ))}{\text{det}(E-(A_B-E)t )}.
\end{eqnarray*}
So $E$ has series Euler characteristic if and only if $-1$ can be substituted to $$\# P^{-1}(b) \frac{\text{sum}(\text{adj}(E-(A_B-E)t ))}{\text{det}(E-(A_B-E)t )}$$ if and only if $-1$ can be substituted to $$ \frac{\text{sum}(\text{adj}(E-(A_B-E)t ))}{\text{det}(E-(A_B-E)t )}$$ if and only if $B$ has series Euler characteristic. Hence, we prove the first claims. If $E$ has series Euler characteristic, then we have 
\begin{eqnarray*}
\chi_{\sum}(E)&=&\# P^{-1}(b) \chi_{\sum}(B) \\
&=& \chi_{\sum} (P^{-1}(b)) \chi_{\sum}(B).
\end{eqnarray*}

The second proof is a proof by the zeta function of a finite category. Suppose 
$$\zeta_B(z)=\prod^n_{k=1} \frac{1}{(1-a_k z)^{b_{k,0}}} \exp\bigg( Q(z)+\sum_{k=1}^n \sum_{j=1}^{e_k-1} \frac{b_{k,j} z^j}{j(1-a_k z)^j} \bigg) $$
for some complex numbers $a_k$ and $b_{k,j}$ and natural numbers $n$ and $e_k$ and a polynomial $Q(z)$ with $\mathbb{Z}$-coefficients whose constant term is 0 (Theorem 3.1 of \cite{NogC}). Then, the uniqueness of the analytic continuation and Theorem \ref{covering} imply
\begin{multline}
\zeta_E(z)=\prod^n_{k=1} \frac{1}{(1-a_k z)^{\# P^{-1}(b) b_{k,0}}} \\ \times\exp\bigg( \# P^{-1}(b)  Q(z)+ \# P^{-1}(b) \sum_{k=1}^n \sum_{j=1}^{e_k-1} \frac{b_{k,j} z^j}{j(1-a_k z)^j} \bigg) .
\end{multline}
We  have $Q(z)=0$ if and only if $\# P^{-1}(b)  Q(z)=0$, so that by Lemma \ref{zeta} the first claim is proven. If $B$ has series Euler characteristic, then we have 
$$\zeta_B(z)=\prod^n_{k=1} \frac{1}{(1-a_k z)^{b_{k,0}}} \exp\bigg( \sum_{k=1}^n \sum_{j=1}^{e_k-1} \frac{b_{k,j} z^j}{j(1-a_k z)^j} \bigg) $$
$$\zeta_E(z)=\prod^n_{k=1} \frac{1}{(1-a_k z)^{\# P^{-1}(b) b_{k,0}}} \exp\bigg( \# P^{-1}(b)  \sum_{k=1}^n \sum_{j=1}^{e_k-1} \frac{b_{k,j} z^j}{j(1-a_k z)^j} \bigg) .$$
Theorem 3.3 of \cite{NogC} implies
\begin{eqnarray*}
\chi_{\sum}(E)&=&\sum^n_{k=1}\sum^{e_k-1}_{j=0} (-1)^j \frac{ \# P^{-1}(b) b_{k,j}}{a_k^{j+1}} \\
&=& \# P^{-1}(b) \sum^n_{k=1}\sum^{e_k-1}_{j=0} (-1)^j \frac{ b_{k,j}}{a_k^{j+1}} \\
&=& \chi_{\sum}(P^{-1}(b)) \chi_{\sum}(B).
\end{eqnarray*}
\end{proof}
\end{them}

\begin{lemma}\label{zeta}
Suppose $I$ is a finite category and the zeta function of $I$ is 
\begin{eqnarray}
\zeta_I(z)=\prod^n_{k=1} \frac{1}{(1-a_k z)^{b_{k,0}}} \exp\bigg( Q(z)+\sum_{k=1}^n \sum_{j=1}^{e_k-1} \frac{b_{k,j} z^j}{j(1-a_k z)^j} \bigg) \label{siki3}
\end{eqnarray} 
for some complex numbers $a_k$ and $b_{k,j}$ and natural numbers $n$ and $e_k$ and a polynomial $Q(z)$ with $\mathbb{Z}$-coefficients whose constant term is 0. Then, $I$ has series Euler characteristic if and only if $Q(z)=0$.
\begin{proof}
Theorem 3.1 of \cite{NogA} implies the zeta function of $I$ is of the form of \eqref{siki3}.

If $I$ has series Euler characteristic, Theorem 3.3 of \cite{NogC} implies $Q(z)=0$. Conversely, let $Q(z)=0$. Here, let us recall what $Q(z)$ is (Theorem 3.1 of \cite{NogC}). Indeed, $Q(z)=\int q(z) dz $ and $q(z)$ is a polynomial with coefficients in $\mathbb{Z}$, moreover
$$\mathrm{sum}(\mathrm{adj}(E-A_I z)A_I)=q(z)|E-A_I z|+r(z)$$
where $$\deg(r(z))<\deg |E-A_Iz|.$$We have $\int q(z) dz =0$ implies $q(z)=0$. Hence, 
$$\mathrm{sum}(\mathrm{adj}(E-A_I z)A_I)=r(z).$$
Lemma 2.3 of \cite{NogC} implies $I$ has series Euler characteristic.
\end{proof}
\end{lemma}

We recall the \textit{Euler characteristic of an $\mathbb{N}$-filtered acyclic category} \cite{Nog11}.

\begin{defi}
A small category $A$ is \textit{acyclic} if every endomorphism and every isomorphism is an identity morphism. 
\end{defi}

\begin{rema}\label{scwol}
This is the same as a skeletal scwol \cite{BH99}. 

Define an order on the set $\mathrm{Ob}(A)$ of objects of $A$ by $x\le y$ if there exists a morphism $x\rightarrow y$. Then, $\mathrm{Ob}(A)$ is a poset.
\end{rema}

\begin{defi}\label{n-fil}
Let $A$ be an acyclic category. A functor $\map{\mu}{A}{\mathbb{N}}$ satisfying $\mu (x)<\mu (y)$ for $x<y$ in Ob$(A)$ is called \textit{an $\mathbb{N}$-filtration of $A$}. A pair $(A, \mu )$ is called \textit{an $\mathbb{N}$-filtered acyclic category}.
\end{defi}

\begin{defi}
Let $(A, \mu)$ be an $\mathbb{N}$-filtered acyclic category. Then, define $\chi_{\text{fil}}(A, \mu)$ as follows. 

For natural numbers $i$ and $n$, let
$$\overline{N_n}(A)_i=\{ \mathbf{f}\in \overline{N_n}(A)\mid \mu(t(\mathbf{f})) =i  \}$$
where $t(\mathbf{f})=x_n$ if $$ \mathbf{f}=\xymatrix{(x_0\ar[r]^{f_1}&x_1\ar[r]^{f_2}&\dots\ar[r]^{f_n}&x_n)}.$$ Suppose each $\overline{N_n}(A)_i$ is finite. Define the formal power series $f_{\chi}(A, \mu)(t)$ over $\mathbb{Z}$ by

$$f_{\chi}(A, \mu)(t)=\sum^\infty_{i=0}(-1)^i\left(\sum^i_{n=0}(-1)^n\# \overline{N_n}(A)_i \right)t^i.$$
Then, define
$$\chi_{\text{fil}}(A, \mu)=f_{\chi}(A, \mu)|_{t=-1}$$
if $f_{\chi}(A,\mu)(t)$ is rational and has a non-vanishing denominator at $t=-1$. 
\end{defi}

\begin{them}
Suppose $(A,\mu_A)$ and $(B,\mu_B)$ are $\mathbb{N}$-filtered acyclic categories and $b$ is  an object of $B$ and $\map{P}{A}{B}$ be a covering whose fiber is finite satisfying the equation $\mu_A(x)=\mu_B(P(x))$ for any object $x$ of $A$. Then, $(A,\mu_A)$ has Euler characteristic $\chi_{\mathrm{fil}}(A,\mu_A)$ if and only if $B$ has Euler characteristic $\chi_{\mathrm{fil}}(B,\mu_B)$. In this case, we have
$$\chi_{\mathrm{fil}}(A,\mu_A)=\chi_{\mathrm{fil}}(P^{-1}(b),\mu) \chi_{\mathrm{fil}}(B,\mu_B)$$
for any $\mathbb{N}$-filtration $\mu$ of $P^{-1}(b)$.
\begin{proof}
For any object $b$ of $B$ and $x$ of $P^{-1}(b)$, we have $$\mu_A^{-1}(i)=\coprod_{b\in \mu_B^{-1}(i)} P^{-1}(b)$$ 
for any $i \ge 0$. Indeed, for any $x$ of $\mu_A^{-1}(i)$, $\mu_B(P(x))=\mu_A(x)=i$. Hence, $P(x)$ belongs to $\mu^{-1}_B(i)$. Moreover, $x$ belongs to $P^{-1}(P(x))$, so that $x$ belongs to $\coprod_{b\in \mu_B^{-1}(i)} P^{-1}(b)$. Conversely, for any $b$ of $\mu_B^{-1}(i)$ and $y$ of $P^{-1}(b)$, $y$ belongs to $\mu_A^{-1}(i)$ since 
$$\mu_A(y)=\mu_B(P(y))=\mu_B(b)=i.$$
Hence, $y$ belongs to $\mu^{-1}_A(i)$. Proposition \ref{n-target} implies
\begin{eqnarray*}
\overline{N_n}(A)_i&=&\coprod_{x\in \mu_A^{-1} (i)} \overline{N_n}(A)_x \\
&=&\coprod_{x\in \coprod_{b\in \mu_B^{-1}(i)} P^{-1}(b)} \overline{N_n}(A)_x \\
&\cong&\coprod_{b\in \mu_B^{-1}(i)} \coprod_{x\in P^{-1}(b)} \overline{N_n}(A)_x \\
&\cong&\coprod_{b\in \mu_B^{-1}(i)} P^{-1}(b) \times \overline{N_n}(B)_b \\
&\cong&P^{-1}(b) \times \coprod_{b\in \mu_B^{-1}(i)} \overline{N_n}(B)_b \\
&=&P^{-1}(b) \times \overline{N_n}(B)_i.
\end{eqnarray*}
Hence, we have
\begin{eqnarray*}
f_{\chi}(A,\mu_A)(t)&=&\sum^\infty_{i=0}(-1)^i\left(\sum^i_{n=0}(-1)^n\# \overline{N_n}(A)_i \right)t^i \\
&=&\sum^\infty_{i=0}(-1)^i\left(\sum^i_{n=0}(-1)^n \# P^{-1}(b) \# \overline{N_n}(B)_i \right)t^i \\
&=&\# P^{-1}(b)f_{\chi}(B,\mu_B)(t).
\end{eqnarray*}
Hence, there exists $\chi_{\mathrm{fil}}(A,\mu_A)$ if and only if the power series $f_{\chi}(A,\mu_A)(t)$ is rational and $-1$ can be substituted to the rational function if and only if the power series $f_{\chi}(B,\mu_B)(t)$ is rational and $-1$ can be substituted to the rational function  if and only if there exists $\chi_{\mathrm{fil}}(B,\mu_B)$. So the first claim is proven. If there exists $\chi_{\mathrm{fil}}(A,\mu_A)$, then we have 
\begin{eqnarray*}
\chi_{\mathrm{fil}}(A,\mu_A)&=&\# P^{-1}(b)\chi_{\mathrm{fil}}(B,\mu_B) \\
&=&\chi_{\mathrm{fil}}(P^{-1}(b),\mu)\chi_{\mathrm{fil}}(B,\mu_B).
\end{eqnarray*}
It is clear that for any $\mathbb{N}$-filtration $\mu$, $\chi_{\mathrm{fil}}(P^{-1}(b),\mu)=\# P^{-1}(b)$. We can give a filtration to $P^{-1}(b)$ , for example, define $\map{\mu}{P^{-1}(b)}{\mathbb{N}}$ by $\mu(x)=0$ for any $x$ of $P^{-1}(b)$.

\end{proof}
\end{them}

\section{Examples}\label{example}

In this section, we give three examples of coverings of small categories.

\begin{exam}\label{exam1}
Let $$\Gamma=\xymatrix{x\ar@<1ex>[r]^f&y\ar@<1ex>[l]^{f^{-1}}}$$ and $B=\mathbb{Z}_2=\{1,-1\}.$ A group can be regarded as a category whose object is just one object $*$ and morphisms are elements of $G$ and composition is the operation of $G$. Define $\map{P}{\Gamma}{B}$ by $P(f)=P(f^{-1})=-1$. Then, $P$ is a covering and this covering was studied in Example 5.33 of \cite{FLS}. Since $\Gamma$ and $B$ are finite groupoids, Proposition 2.3 of \cite{NogA} implies
$$\zeta_{\Gamma}(z)=\frac{1}{(1-2z)^2}, \zeta_{B}(z)=\frac{1}{1-2z}.$$
The number of sheet of $P$ is 2. We have $\zeta_{\Gamma}(z)=\zeta_B(z)^2$. 

Applying the classifying space functor $B$ to $P$, we obtain the famous covering $\map{p_{\infty}}{\mathbb{S}^{\infty}}{\mathbb{RP}^{\infty}}$

$$\xymatrix{B\Gamma \ar[r]^{\cong} \ar[d]^{BP}&\mathbb{S}^{\infty}\ar[d]^{p_{\infty}} \\ 
B\mathbb{Z}_2\ar[r]^{\cong}&\mathbb{RP}^{\infty}.}$$
For a finite groupoid $G$, three Euler characteristic $\chi_L,\chi_{\sum}$ and $\chi^{(2)}$ coincide and the value is $$\sum_{x\in \mathrm{Ob}(G)/\cong} \frac{1}{\#\mathrm{Aut}(x)}$$
where this sum runs over all the isomorphism classes of objects of $G$ (Example 2.7 of \cite{Leia}, Theorem 3.2 of \cite{Leib} and Example 5.12 of \cite{FLS}). Hence, we have 
$$\chi_{\sum}(\Gamma)=1, \chi_{\sum}(B)=\frac{1}{2}, \chi_{\sum}(P^{-1}(*))=2$$
and $$\chi_{\sum}(\Gamma)=\chi_{\sum}(P^{-1}(*))\chi_{\sum}(B).$$
\end{exam}

\begin{exam}
Let $$A=\xymatrix{y_1&&y_2&&y_3&\dots&y_n \\
&x_1\ar[ul]^{f_1}\ar[ur]_{g_1}&&x_2\ar[ul]^{f_2}\ar[ur]_{g_2}&\dots&x_{n-1}\ar[ul] \ar[ur]^{g_{n-1}}&&x_n\ar[ulllllll]_(0.3){g_n} \ar[ul]^{f_n}}$$
and $$B=\xymatrix{a\ar@/^/[r]^{h_1} \ar@/_/[r]_{h_2}&b}.$$ Define a functor $\map{P}{A}{B}$ by $P(x_i)=a$,  $P(y_i)=b$, $P(f_i)=h_1$ and $P(g_i)=h_2$ for any $i$. Then, $P$ is a covering. By Proposition 2.9 of \cite{NogA}, we have 
$$\zeta_A(z)=\frac{1}{(1-z)^{2n}} \exp\left(\frac{2n z}{1-z}\right), \zeta_B(z)=\frac{1}{(1-z)^{2}} \exp\left(\frac{2 z}{1-z}\right).$$
The number of sheet of $P$ is $n$. We have $\zeta_A(z)=\zeta_B(z)^n$.

Applying the classifying space functor $B$ to $P$, we obtain the famous covering $\map{p^n}{\mathbb{S}^1}{\mathbb{S}^1}$
where $$\mathbb{S}^1=\{z\in \mathbb{C} \mid |z|=1\}$$ and $p^n$ is the $n$-th power mapping. The map $p^n$ is a covering (see \cite{Hat02}, for example)

$$\xymatrix{B A\ar[r]^{\cong} \ar[d]^{BP}&\mathbb{S}^{1}\ar[d]^{p^n} \\ 
B B\ar[r]^{\cong}&\mathbb{S}^1.}$$

The two categories $A$ and $B$ are finite acyclic categories. For a finite acyclic category, four Euler characteristics $\chi_L,\chi_{\sum},\chi^{(2)}$ and $\chi_{\mathrm{fil}}$ coincide (see, for example, Introduction of \cite{Nog}). Furthermore, they coincide the Euler characteristic for cell complexes of the classifying space of an acyclic category (Proposition 2.11 of \cite{Leia}). We have 
$$\chi_{\sum}(A)=0, \chi_{\sum}(B)=0, \chi_{\sum}(P^{-1}(a))=2$$
and 
$$\chi_{\sum}(A)=\chi_{\sum}(P^{-1}(a))\chi_{\sum}(B).$$
\end{exam}

We introduce an example of a covering of infinite categories.

\begin{exam}
Suppose
$$A=\xymatrix{x_0\ar[r]\ar[dr]&x_1\ar[r]\ar[dr]&x_2\ar[r]\ar[dr]&\cdots \\ y_0\ar[r]\ar[ur]&y_1\ar[r]\ar[ur]&y_2\ar[r]\ar[ur]&\cdots }$$
and 
$$B=\xymatrix{b_0\ar@<1ex>[r]\ar@<-1ex>[r]&b_1\ar@<1ex>[r]\ar@<-1ex>[r]&b_2\ar@<1ex>[r]\ar@<-1ex>[r]&\cdots}$$
Where $A$ is a poset, that is, $A$ acyclic and each hom-set has at most exactly one morphism and for $n<m$ and $b_n$ and $b_m$, define
$$\mathrm{Hom}_B(b_n,b_m)=\{ \varphi^0_{n,m}, \varphi^1_{n,m} \}$$
and a composition of $B$ is defined by $\varphi^j_{m,\ell}\circ \varphi^i_{n,m}=\varphi^k_{n,\ell }$ where $k \equiv i+j \mod 2$ for $n<m<\ell$. Define $\map{P}{A}{B}$ by $P(x_i)=P(y_i)=b_i$ and $P((x_n,x_m))=P((y_n,y_m))=\varphi^0_{n,m}$ and $P((y_n,x_m))=P((x_n,y_m))=\varphi^1_{n,m}$ for $n<m$. Then, $P$ is a covering. The indexes of objects of $A$ and $B$ give $\mathbb{N}$-filtrations $\mu_A$ and $\mu_B$ to $A$ and $B$, respectively. We have
\begin{eqnarray*}
f_{\chi}(A,\mu_A)(t)&=&\sum^{\infty}_{i=0} \bigg( \sum^i_{n=0} (-1)^n 2^{n+1}\binom{i}{n}\bigg) t^i \\
&=&2\sum^{\infty}_{i=0} t^i \\
&=&\frac{2}{1-t},
\end{eqnarray*}
so that $\chi_{\mathrm{fil}}(A,\mu_A)=1.$
We have 
\begin{eqnarray*}
f_{\chi}(B,\mu_B)(t)&=&\sum^{\infty}_{i=0} \bigg( \sum^i_{n=0} (-1)^n 2^{n}\binom{i}{n}\bigg) t^i \\
&=&\sum^{\infty}_{i=0} t^i \\
&=&\frac{1}{1-t},
\end{eqnarray*}
so that $\chi_{\mathrm{fil}}(B,\mu_B)=\frac{1}{2}.$
In fact, the categories $A$ and $B$ are the barycentric subdivision of $\Gamma$ of Example \ref{exam1} and $\mathbb{Z}_2$ (see \cite{Nog11} and \cite{Nog}). Hence, Theorem 4.9 of \cite{Nog11} and Example \ref{exam1} imply their Euler characteristic $\chi_{\mathrm{fil}}(A,\mu_A)$ and $\chi_{\mathrm{fil}}(B,\mu_B)$. We obtain 
$$\chi_{\mathrm{fil}}(A,\mu_A)=\chi_{\mathrm{fil}}(P^{-1}(b_0),\mu)\chi_{\mathrm{fil}}(B,\mu_B)$$
for any $\mathbb{N}$-filtration $\mu$ of $P^{-1}(b)$.

\end{exam}

\end{document}